\documentclass[10pt]{amsart}
\usepackage{amsmath,mathrsfs,amssymb}
\usepackage{fullpage}
\newtheorem{teo}{Theorem}
\newtheorem{pro}{Proposition}
\newtheorem{lem}{Lemma}
\newtheorem{cor}{Corollary}
\newtheorem*{rem}{Remark}

\title{Moments of the Hermitian Matrix Jacobi process} 
\author[L. Deleaval]{Luc Deleaval}
\address{LAMA, Universit\'e Marne la Vall\'ee \\ Champs sur Marne \\ 77454 Marne la VallŽe Cedex 2, France}
\email{luc.deleaval@u-pem.fr}
\author[N. Demni]{Nizar Demni}
\address{IRMAR, Universit\'e de Rennes 1\\ Campus de
Beaulieu\\ 35042 Rennes cedex\\ France}
\email{nizar.demni@univ-rennes1.fr}


\begin{document}
\maketitle

\begin{abstract}
In this paper, we compute the expectation of traces of powers of the hermitian matrix Jacobi process for a large enough but fixed size. To proceed, we first derive the semi-group density of its eigenvalues process as a bilinear series of symmetric Jacobi polynomials. Next, we use the expansion of power sums in the Schur polynomial basis and the integral Cauchy-Binet formula in order to determine the partitions having non zero contributions after integration. It turns out that these are hooks of bounded weight  and the sought expectation results from the integral of a product of two Schur functions with respect to a generalized Beta distribution. For special values of the parameters on which the matrix Jacobi process depends, the last integral reduces to the Cauchy determinant and we close the paper with the investigation of the asymptotic behavior of the resulting formula as the matrix size tends to infinity. 
\end{abstract}

\section{Reminder and motivation} 
Given three integers $d,p,m$ such that $1 \leq p,m < d$, the hermitian matrix Jacobi process $(J_t)_{t \geq 0}$ of parameters $(p,q:=d-p)$ was defined in \cite[page 141]{Dou} as the product of the $m \times p$ upper-left corner of a $d \times d$ Brownian motion $(Y_t)_{t \geq 0}$ on the unitary group $\mathcal{U}(d, \mathbb{C})$ (\cite{Lia}) and of its Hermitian conjugate. Equivalently, if $P_m$ and $Q_p$ are two $d \times d$ diagonal projections of ranks $m$ and $p$ respectively, then
\begin{equation*}
\left(\begin{array}{cc}
J_t &  0_{m,d-m} \\
0_{d-m,m} & 0_{d-m}\end{array}\right) :=(P_mY_tQ_p)(P_mY_tQ_p)^{\star} = P_mY_tQ_pY_t^{\star}P_m,
\end{equation*}
where $0_{d-m,m}, 0_{m,d-m}, 0_{d-m}$ are the null matrices of shapes $d-m \times m, m \times d-m,$ and $d-m\times d-m$ respectively. With this matrix representation in hands and from the independence of the increments of the L\'evy process $(Y_t)_{t \geq 0}$, it follows that if $d=d(m)$ and $p=p(m)$ depend on $m$ such that
\begin{equation}\label{Limits}
\lim_{m \rightarrow \infty} \frac{p(m)}{d(m)} := \theta \in ]0,1[, \quad \lim_{m \rightarrow \infty} \frac{m}{p(m)} := \eta >0,  \,\,\textrm{with} \, \, \eta \theta \in ]0,1[,
\end{equation}
exist, then the expectation of the normalized trace of any finite-tuple of matrices drawn from $(J_{t/d(m)})_{t \geq 0}$ converge as $m \rightarrow \infty$ (\cite{Dem00}, see also the recent paper \cite{DCK} where the convergence is shown to hold in the strong sense). In particular, if $\mathbb{E}$ denotes the expectation of the probability space where $(Y_t)_{t \geq 0}$ is defined, then the following limit
\begin{equation}\label{limit} 
\mathcal M_n(t, \eta, \theta):= \lim_{m \rightarrow \infty} \frac{1}{m}\mathbb{E}\left(\mathrm{tr}\left[\left(J_{t/d(m)}\right)^n\right]\right)
\end{equation}
exists for any $n \geq 0, \, t \geq 0$, and the sequence $(\mathcal M_n(t,\eta,\theta))_{n \geq 0}$ determines the spectral distribution of the so-called free Jacobi process (\cite{Dem00}). Furthermore, the limit
\begin{equation*}
\mathcal M_n(\infty, \eta, \theta):=\lim_{t \rightarrow \infty} \mathcal M_n(t, \eta, \theta), \quad n \geq 0,
\end{equation*}
is the moment sequence of the spectral distribution of the large $m$-limit of $P_mUQ_{p(m)}U^{\star}P_m$, where $U$ is a $d(m) \times d(m)$ Haar unitary matrix. In other words, this limiting distribution describes the spectrum of the large $m$-limit of matrices drawn from the Jacobi unitary ensemble and its Lebesgue decomposition follows readily from freeness considerations (\cite{Cap-Cas}, \cite{Col}, \cite{Dem00}). Besides, an explicit expression of 
$\mathcal{M}_n(\infty, \eta, \theta)$ obtained from large $m$-asymptotics of the moments of the multivariate Beta distribution figures in \cite[Theorem 4.4]{CDLV}. However, the situation becomes rather considerably more complicated when dealing with 
$\mathcal M_n(t, \eta, \theta)$ for fixed time $t > 0$, as witnessed by the series of papers \cite{DHH}, \cite{DH} and \cite{Dem-Ham}. For instance, it was proved in \cite{DHH} that
\begin{equation}\label{Mom}
\mathcal M_n(t,1,1/2) =  \frac{1}{2^{2n}}  \binom{2n}{n} + \frac{1}{2^{2n-1}} \sum_{k=1}^n  \binom{2n}{n-k}\frac{1}{k} L_{k-1}^{1}(2kt) e^{-kt},
\end{equation}
where $L_k^1$ is the $k$-th Laguerre polynomial of index $1$ (\cite[chapter 6]{AAR}). In this formula
\begin{equation*}
\frac{1}{k} L_{k-1}^{1}(2kt) e^{-kt}, \quad k \geq 1,
\end{equation*}
is the $k$-th moment of the so-called free unitary Brownian motion at time $2t$ (\cite{Biane}, \cite{Levy}, \cite{Rains}), which arises in the large $d$-limit of $(Y_{t/d})_{t \geq 0}$. This observation led to a beautiful, yet striking, representation of the spectral distribution of the free Jacobi process associated with the couple of values $\eta = 1, \theta = 1/2$. In \cite{DH}, partial results on the spectrum of the free Jacobi process associated with $\eta = 1,\theta \in ]0,1]$ were obtained. There, a unitary process related to the free Jacobi process was considered and a detailed analysis of the dynamics of its spectrum was performed. The connection between both spectra is then ensured by a non commutative binomial-type expansion. In the recent paper \cite{Dem-Ham}, a complicated expression of $\mathcal M_n(t, 1, \theta)$ is obtained using sophisticated tools from complex analysis. 

Motivated by these findings, we tackle here the problem of computing the large $m$-limit \eqref{limit} by deriving an explicit expression of $\mathbb{E}(\mathrm{tr}[(J_{t/d})^n])$ for fixed $t > 0, n \geq 1$. To this end, we shall assume that $m$ is large enough so that $m > n$ and make use of the semi-group density of the eigenvalues process of $(J_t)_{t \geq 0}$. In this respect, it was noticed in \cite{Dou} that the latter process is realized as $m$ independent real Jacobi processes of parameters $(2(p-m+1) > 0 , 2(q-m+1) > 0)$ and conditioned never to collide. As a matter of fact, its semi-group density follows readily from the Karlin and McGregor formula (see \cite{Dem01} for the details) and is given by a bilinear series of symmetric Jacobi polynomials indexed by partitions (see for instance \cite{Deb}, \cite{Las1}). We shall also prove the absolute-convergence of the series defining this density, so that Fubini Theorem applies when computing $\mathbb{E}(\mathrm{tr}[(J_{t/d})^n])$. Next, with the help of the expansion of the $n$-th power sum 
in the Schur polynomial basis (\cite{MD}) and of the integral Cauchy-Binet formula (\cite[page 37]{DG}), we determine the partitions having non zero contributions after integration. These are exactly the hooks of weights less than $n$, and both papers \cite{Las0} and \cite{Las1} provide an explicit expansion of the corresponding symmetric Jacobi polynomial in the Schur polynomial basis. The sought expectation follows from the integral of a product of Schur functions with respect to a multivariate Beta weight. The Cauchy-Binet formula allows once more to express this integral as a determinant of a matrix whose entries are Beta functions (see Exercise 8, p.386 in \cite{MD}). Summarizing, we obtain the following result, where we denote by $U_{\tau}^{p-m,q-m,m}(1^m)$ the value at the point 
\begin{equation*}
1^m = (\underbrace{1,\ldots,1}_{m \mathrm{\,times}})
\end{equation*}
 of the symmetric Jacobi polynomial $U_{\tau}^{p-m,q-m,m}$ of parameters $(p-m, q-m)$, by $\mu \subseteq \tau$ the ordering induced by the Young diagrams associated with the partitions $\mu=(\mu_1\geq\ldots\geq\mu_m\geq0)$, $\tau=(\tau_1\geq\ldots\geq\tau_m\geq0)$, by  $\alpha = \alpha(n,k) := (n-k, 1^k)$ a hook of weight $|\alpha| = n$ and by $\beta(\cdot,\cdot)$  the Beta function (see the next sections for more details on both Jacobi polynomials and partitions).
\begin{teo}\label{Obs}
Let $p \wedge q \geq m$ and set $r:=p-m \geq 0$, $s:=q-m \geq 0$. Then 
\begin{multline}\label{princi}
\mathbb{E}(\mathrm{tr}[(J_{t/d})^n])=\\ \sum_{k=0}^{n-1}(-1)^k\sum_{\tau \subseteq \alpha}a_{\tau}^{r,s,m} e^{-K_{\tau}^{r,s,m}(t/d)}U_{\tau}^{r,s,m}(1^m) 
\sum_{\mu \subseteq \tau}b_{\mu,\tau}^{r,s,m}  \det\left(\beta(\alpha_i + \mu_j+2m - i- j +r +1,s+1)\right)_{i,j=1}^m,
\end{multline}
where $a_{\tau}^{r,s,m}, b_{\mu,\tau}^{r,s,m}  \in \mathbb{R}$ are given in \eqref{coeff1} and \eqref{coeff2} respectively and where 
\[
K_{\tau}^{r,s,m}=\sum_{i=1}^m\tau_i(\tau_i + r+s+1 + 2(m-i)).
\]  
\end{teo}

When $s=0$, the determinant of Beta functions reduces to the well-known Cauchy determinant. Together with Weyl dimension formula, we get the following corollary where, for a partition $\tau$, $s_{\tau}$ denotes the associated Schur polynomial (see Section 3 for more details on Schur polynomial).
\begin{cor}\label{coll}
If $s=0$, then we have
\begin{align}\label{Formula1}
\mathbb{E}(\mathrm{tr}[(J_{t/d})^n]) &= \sum_{k=0}^{n-1}(-1)^k\sum_{\tau \subseteq \alpha}e^{-K_{\tau}^{r,0,m}(t/d)} \prod_{i=1}^m[2(\tau_i+m-i)+r+1]\left\{\frac{\Gamma(m-i+1)\Gamma(r+\tau_i+m-i+1)}{\Gamma(\tau_i+m-i+1)\Gamma(r+m-i+1)}\right\}^2 
\nonumber
\\& [s_{\tau}(1^m)]^2 U_{\tau}^{r,0,m}(1^m) 
 \sum_{\mu \subseteq \tau}  \frac{\prod_{1 \leq i < j \leq m}(\tau_i+\tau_j+2m-i-j+r+1)^2}{\prod_{i,j=1}^m(\alpha_i+\mu_j+2m-i-j+r+1)} b_{\mu,\tau}^{r,0,m}s_{\mu}(1^m) s_{\alpha}(1^m).
\end{align}
If further $r=s=0$, then
\begin{align*}
\mathbb{E}(\mathrm{tr}[(J_{t/d})^n]) &= 
\sum_{k=0}^{n-1}(-1)^k\sum_{\tau \subseteq \alpha}e^{-K_{\tau}^{0,0,m}(t/d)} \prod_{i=1}^m[2(\tau_i+m-i)+1]
\\& [s_{\tau}(1^m)]^2 U_{\tau}^{0,0,m}(1^m) 
 \sum_{\mu \subseteq \tau} \frac{\prod_{1 \leq i < j \leq m}(\tau_i+\tau_j+2m-i-j+1)^2}{\prod_{i,j=1}^m(\alpha_i+\mu_j+2m-i-j+1)} b_{\mu,\tau}^{0,0,m} s_{\mu}(1^m) s_{\alpha}(1^m).
\end{align*}
\end{cor}
Let us point out that, for $s=1$, the determinant 
\begin{multline*}
\det\left(\beta(\alpha_i + \mu_j+2m - i- j +r +1,2)\right)_{i,j=1}^m\\=
\det\left(\frac{1}{(\alpha_i + \mu_j+2m - i- j +r +1)(\alpha_i + \mu_j+2m - i- j +r +2)}\right)_{i,j=1}^m
\end{multline*}
was already considered in \cite{Lasc}, where it is expanded in some basis of symmetric functions. Up to our best knowledge, there is no general explicit expression of the above determinant for arbitrary $s \geq 0$. Nonetheless, as we shall see below, the term corresponding to the null partition $\tau= (0)$ may be computed using Kadell's integral (see Exercise 7, p.385 in \cite{MD}) and as such, we retrieve the moments derived in \cite{CDLV} (see Proposition 2.2 and Corollary 2.3 there) of the multivariate Beta distribution arising from the Jacobi unitary ensemble. This is by no means a surprise since all but this term cancel when we let $t \rightarrow \infty$ in \eqref{princi} and the distribution of $J_t$ converges weakly as $t \rightarrow \infty$ to that of the Jacobi unitary ensemble (also known as the matrix-variate Beta distribution).   

Back to formula \eqref{Formula1}, some of the products involved there terminate after cancellations, since the lengths of $\mu, \tau, \alpha$ satisfy $l(\mu) \leq l(\tau) \leq l(\alpha) \leq n< m$. This observation allows to take the limit as $m \rightarrow \infty$ there, assuming $p = p(m)$ and $d = d(m)$ are such that \eqref{Limits} holds. Moreover, we show that $b_{\mu,\tau}^{r(m),s(m),m}s_{\mu}(1^m)$ has finite large $m$-limit which, together with the generalized binomial formula for Schur functions (\cite{Lass}), entail
\begin{equation*}
\lim_{m \rightarrow \infty}U_{\tau}^{r(m),s(m),m}(1^m) = \left(1-\frac{1}{\theta}\right)^{|\tau|}.
\end{equation*}
Here, we write $r=r(m) = p(m) - m, s = s(m) = d(m)-p(m) -m$ and the assumption $s(m) =  0$ corresponds in the large $m$-limit to the set 
\begin{equation*}
\{\theta \in ]0,1[, \theta(1+\eta) = 1\}. 
\end{equation*}
Since $s_{\alpha}(1^m) = O(m^{|\alpha|})$ for any partition $\alpha$ (\cite{CDLV}, p.4), then we are led after normalizing by the factor $(1/m)$ to an indeterminate limit and as such, the computation of \eqref{limit} seems to be out of reach for the moment. Note in this respect that the derivation of the moments $\mathcal M_n(\infty, \eta, \theta)$ performed in \cite{CDLV} is based on the inverse binomial transform. 

The paper is organized as follows. In the next section, we recall the definition of the Brownian motion on the unitary group $\mathcal{U}(d, \mathbb{C})$ and derive the stochastic differential equation satisfied by the hermitian matrix Jacobi process which was announced in \cite{Dou} without proof. In the same section, we recall also the stochastic differential system satisfied by the corresponding eigenvalues process and prove the absolute convergence of the semi-group density of the latter. In section 3, we prove our main results, that is,Theorem \ref{Obs} and his corollary. For that purpose, we recall some facts on both Schur polynomials and symmetric Jacobi polynomials associated with hooks then generalize an orthogonality relation for the real Jacobi polynomial to its multivariate analogue. In the last section, we investigate the asymptotic behavior of all the terms appearing in the right-hand side of \eqref{Formula1}.

\section{The Hermitian matrix Jacobi process and its eigenvalues process}
\subsection{From the unitary Brownian motion to the Hermitian matrix Jacobi process} 
The existence of the limit \eqref{limit} relies to a large extent on the convergence of the moments of $(Y_{t/d})_{d \geq 0}$ to those of the free unitary Brownian motion (\cite{Biane}). The time normalization $t/d$ is equivalent to the normalization of the Laplace-Beltrami operator on $\mathcal{U}(d)$ by a factor $1/d$, which in this case corresponds to the Killing form
\begin{equation*}
-d\,\mathrm{tr}(XY),
\end{equation*}
where $X,Y$ are skew-hermitian matrices. With this normalization, the unitary Brownian motion solves the following stochastic differential equation (see \cite{Lia}): 
\begin{equation}\label{UBM}
dY_t = iY_tdH_t - \frac{1}{2}Y_t dt, \quad Y_0 = {\it I}_d,
\end{equation}
where ${\it I}_d$ is the $d \times d$ identity matrix and $(H_t)_{t \geq 0}$ is a $d \times d$ matrix-valued Hermitian process whose diagonal entries are real Brownian motions while its off diagonal entries are complex Brownian motions, all of them being independent and have common variance $t/d$. Besides, the process $(Y_{t/d})_{t \geq 0}$ is a left Brownian motion in the sense that the semi-group operator 
\begin{equation*} 
f \mapsto \mathbb{E}\left[f(ZY_t)\right], \qquad Z \in \mathcal{U}(d,\mathbb{C}), 
\end{equation*}
defined on the space of continuous functions $f$ on $\mathcal{U}(d,\mathbb{C})$ is left-invariant. Equivalently, the right-increments  
\begin{equation*}
Y_{s/d}^{-1}Y_{t/d}, \quad 0 \leq s < t, 
\end{equation*}
of $(Y_{t/d})_{t \geq 0}$ are invariant under left multiplication by any complex unitary matrix. This choice is by no means a loss of generality since the process $(Y_{t/d}^{-1})_{t \geq 0}$ has the same distribution as $(Y_{t/d})_{t \geq 0}$ and is a right Brownian motion on $\mathcal{U}(d,\mathbb{C})$. 

Now, one can use \label{UBM} in order to derive a stochastic differential equation satisfied by $J_t$. To this end, let
\begin{equation*}
Y_t = \left(\begin{array}{lr}
X_t & U_t \\ 
V_t & W_t \end{array}\right),  \quad H_t = \left(\begin{array}{lr}
R_t & S_t \\ 
M_t & N_t \end{array}\right) ,
\end{equation*}
be the block decompositions of $Y_t$ and $H_t$. Here, $X_t$ is the $m \times p$ upper-left corner of $Y_t$ so that $J_t = X_tX_t^{\star}$, while $U_t, V_t, W_t, R_t, S_t, M_t, N_t$ are $m \times q, d-m \times p, d-m \times q, p \times p, p \times q, q \times p, q \times q$ matrices respectively. Hence, \eqref{UBM} readily gives 
\begin{equation*}
dX_t = i(X_tdR_t + U_tdM_t) - \frac{X_t}{2} dt
\end{equation*}
and It\^o formula yields 
\begin{equation*}
dJ_t  = X_t(dX_t^{\star}) + (dX_t) X_t^{\star} + < (dX_t),(dX_t^{\star})> 
\end{equation*}
where $<\cdot,\cdot>$ denotes the bracket of continuous semi-martingales (\cite{Rev-Yor}). Since
\begin{equation*}
<dB_t, d\overline{B_t}> = t, \quad <dB_t, dB_t> = 0, 
\end{equation*}
for any complex Brownian motion $(B_t)_{t \geq 0}$ of variance $t$, since $(R_t)_{t \geq 0}$ and $(M_t)_{t \geq 0}$ are independent and since $X_tX_t^{\star} + U_tU_t^{\star} = {\it I}_m$, then the finite-variation part of the semi-martingale decomposition of $dJ_t$ is given by:
\begin{equation*}
\left[\frac{p}{d}{\it I}_m - J_t\right]dt. 
\end{equation*}
Again, since $R_t$ is Hermitian, then the local-martingale part of $dJ_t$ is given by 
\begin{equation*}
i\left[U_tdM_tX_t^{\star} - X_tdM_t^{\star} U_t^{\star}\right],  
\end{equation*}
whose bracket coincides with that of the local martingale:
\begin{equation*}
\sqrt{J_t}dF_t\sqrt{{\it I}_m - J_t} + \sqrt{{\it I}_m - J_t}dF_t^{\star}\sqrt{J_t},
\end{equation*}
where $(F_t)_{t \geq 0}$ is a complex Brownian matrix whose entries are independent and have common variance $t/d$. Hence, if $J_0$ and ${\it I}_m - J_0$ are positive-definite, the following stochastic differential equation holds
\begin{equation*}
dJ_t = \sqrt{J_t}dF_t\sqrt{{\it I}_m - J_t} + \sqrt{{\it I}_m - J_t}dF_t^{\star}\sqrt{J_t} + \left(\frac{p}{d}{\it I}_m - J_t\right)dt
\end{equation*}
as long as $J_t$ and ${\it I}_m - J_t$ remain so. According to Bru's Theorem (see \cite[page 3061]{Kat-Tan}), there exist $m$ real Brownian motions $(\nu_i)_{i=1}^m$ with common variance $t$ such that the eigenvalues process, say $(\lambda_i)_{i =1}^m$, satisfies the stochastic differential system 
\begin{align}\label{DiffSys}
d\lambda_i(t) = \sqrt{(2/d)(\lambda_i(t)(1-\lambda_i(t))}d\nu_i(t) + \frac{1}{d}\left[(p - d\,\lambda_i(t)) + 
 \sum_{j \neq i}\frac{\lambda_i(t)(1-\lambda_j(t)) + \lambda_j(t)(1-\lambda_i(t))}{\lambda_i(t) - \lambda_j(t)}\right] dt
\end{align}
as long as $0 < \lambda_m(t) < \lambda_{m-1}(t) < \dots < \lambda_1(t) < 1$. Recalling $q=d-p$, then the infinitesimal generator of $(\lambda_i(2td), t \geq 0)_{i=1}^m$ coincides with the one displayed in \cite[page 150]{Dou}. Consequently, $(\lambda_i)_{i=1}^m$ is realized as a Doob transform of $m$ independent real Jacobi processes of parameters $(2(p-m+1), 2(q-m+1))$ killed when they first collide, the sub-harmonic function being the Vandermonde polynomial. On the other hand, the main result proved in \cite{Dem01} shows that if $p \wedge q > m-(1/2)$, then \eqref{DiffSys} admits, for any starting point $\lambda(0)=(0\leq\lambda_m(0) \leq \dots \leq \lambda_1(0)\leq1)$, a unique strong solution defined on the whole positive half-line. Altogether, we deduce from the last section of \cite{Dem01} that the semi-group density\footnote{With respect to Lebesgue measure $d\lambda = \prod_{i=1}^md\lambda_i$.} of $(\lambda_i)_{i =1}^m$, say $G_t^{r,s,m}$, is given at time $t$ by: 
\begin{align}\label{Jac-Sem}
G_t^{r,s,m}(\lambda(0),\lambda) = \sum_{\tau = (\tau_1 \geq \dots \geq \tau_m \geq 0)}e^{-K_{\tau}^{r,s,m}(t/d)} \frac{\det[P_{\tau_i + m-i}^{r,s}(\lambda_j(0))]_{i,j=1}^m}{V(\lambda(0))}\frac{\det[P_{\tau_i+m-i}^{r,s}(\lambda_j)]_{i,j=1}^m}{V(\lambda)}W^{r,s,m}(\lambda),
\end{align}
where we recall $r = p-m, s = q-m$, 
\[
K_{\tau}^{r,s,m} := \sum_{i=1}^m\tau_i(\tau_i + r+s+1 + 2(m-i)),
\]
where we have set
\begin{eqnarray*}
V(\lambda) &:=& \prod_{i < j}(\lambda_i - \lambda_j), \\ 
W^{r,s,m}(\lambda) &:=& \prod_{i=1}^m\lambda_i^{r}(1-\lambda_i)^{s} V(\lambda)^{2}{\bf 1}_{\{0 < \lambda_m < \dots < \lambda_1 < 1\}}, \\
\end{eqnarray*}
and where $P_n^{r,s}$ stands for the $n$-th orthonormal Jacobi polynomial on $[0,1]$. Actually, 
\begin{equation*}
P_n^{r,s} := \frac{p_{n}^{r,s}}{\|p_{n}^{r,s}\|_2}= \frac{1}{\|p_{n}^{r,s}\|_2}\frac{(r+1)_n}{n!}{}_2F_1(-n,n+r+s+1,r+1,\mathopen\cdot\mathclose)
\end{equation*}
with 
\[
\|p_n^{r,s}\|_2^2 := \frac{1}{2n+r+s+1}\, \frac{\Gamma(r+n+1)\Gamma (s+n+1)}{\Gamma(n+1)\Gamma(n+1+r+s)}, \quad (r+1)_n = \frac{\Gamma(r+1+n)}{\Gamma(r+1)},
\]
and ${}_2F_1$ is the Gauss hypergeometric function (see \cite[chapters 2 and  6]{AAR} for more details). Set
\begin{equation*}
P_{\tau}^{r,s,m}(x):=\frac{\det[P_{\tau_i + m-i}^{r,s}(x_j)]_{i,j=1}^m}{V(x)}=\prod_{i=1}^m \frac{1}{\|p_{\tau_i+m-i}^{r,s}\|_2} \frac{\det[p_{\tau_i + m-i}^{r,s}(x_j)]_{i,j=1}^m}{V(x)},
\end{equation*}
then $P_{\tau}^{r,s,m}$ is known as the symmetric (orthonormal) Jacobi polynomial associated with the partition $\tau$. 
Under different normalizations, the family $(P_{\tau}^{r,s,m})_{\tau}$ appeared independently in \cite{Bee-Opd}, \cite{Deb}, \cite{Las1}, \cite{Ols-Oko} and \cite{Ols-Osi}. 
For instance, since 
\begin{equation*}
p_{\tau_i+m-i}^{r,s}(0)= \frac{(r+1)_{\tau_i+m-i}}{(\tau_i+m-i)!},
\end{equation*}
 then $G_t^{r,s,m}(\lambda(0),\lambda)$ may be written as
\begin{multline}\label{Jac-Sem1}
G_t^{r,s,m}(\lambda(0),\lambda) = \\
\sum_{\tau = (\tau_1 \geq \dots \geq \tau_m \geq 0)}e^{-K_{\tau}^{r,s,m}(t/d)}
\left\{V(\tilde{\tau}) \prod_{1 \leq i < j \leq m} \frac{1}{(r+j-i)i} \prod_{i=1}^m \frac{ p_{\tau_i+m-i}^{r,s}(0)}{||p_{\tau_i+m-i}^{r,s}||_2}\right\}^2U_{\tau}^{r,sm}(\lambda(0))U_{\tau}^{r,s,m}(\lambda)W^{r,s,m}(\lambda) ,
\end{multline}
where $U_{\tau}^{r,s,m}$ denotes the polynomial considered in \cite{Las1}, normalized to be equal to $1$ at $(\underbrace{0,\dots,0}_{m \mathrm{\, times}})$, see \cite[ Theorem 10]{Las1}. More explicitely
\begin{equation*}
U_{\tau}^{r,s,m}(\lambda):=\frac{ (-1)^{m(m-1)/2}}{V(\tilde{\tau})}\prod_{1 \leq i < j \leq m}(r+j-i)i\,\,\frac{\det({}_2F_1(-(\tau_i+m-i),\tau_i+m-i+r+s+1,r+1,\lambda_j))_{i,j=1}^m}{V(\lambda)}
\end{equation*}
with 
\[
V(\tilde{\tau}) = \prod_{1 \leq i < j \leq m}(\tau_i-\tau_j+j-i)(\tau_i+\tau_j+2m-i-j+r+s+1). 
\]
The representation \eqref{Jac-Sem1} is convenient for our purposes since when $\tau$ is a hook, an explicit expansion of $U_{\tau}^{r,s,m}$ in the Schur polynomial basis is given in \cite{Las1}.   

\subsection{Absolute convergence of the semi-group density} 
Another normalization of the symmetric Jacobi polynomial is related to the spherical function property they satisfy for special parameters $(r,s)$ (see Table II in \cite{Ols-Oko}). It has the merit to be well-suited for proving that the series given in \eqref{Jac-Sem} is absolutely convergent. Indeed, let $\phi \in [-1,1]^m$ and let
\begin{equation*}
q_n^{r,s}(x) = p_n^{r,s}((1-x)/2), \quad Q_{\tau}^{r,s,m}(\phi) = \frac{\det[q_{\tau_i + m-i}^{r,s}(\phi_j)]_{i,j=1}^m}{V(\phi)},
\end{equation*}
be the Jacobi polynomial in $[-1,1]$ and the symmetric Jacobi polynomial in $[-1,1]^m$ respectively. Then Proposition 7.2 in \cite{Ols-Oko} shows that $Q_{\tau}^{r,s,m}$ coincides up to a constant with the symmetric Jacobi polynomials considered there. Moreover, Proposition 1.1 in the same paper shows that for any $\phi \in [-1,1]^m$, 
\begin{equation*}
|Q_{\tau}^{r,s,m}(\phi)| \leq Q_{\tau}^{r,s,m}(1^m), \quad  \quad r \geq s \geq 0,
\end{equation*}
while the special value $Q_{\tau}^{r,sm}(1^m)$ is given by (\cite[Proposition 7.1]{Ols-Oko}):
\begin{equation*}
Q_{\tau}^{r,s,m}(1^m) = V(\tilde{\tau}) \prod_{i=1}^m \frac{\Gamma(\tau_i+m-i+r+1)2^{-(m-i)}}{\Gamma(\tau_i+m-i+1)\Gamma(m-i+r+1)\Gamma(m-i+1)}.
\end{equation*}
Since 
\begin{equation*}
P_{\tau}^{r,s,m}(x) = (-2)^{m(m-1)/2}\prod_{i=1}^m \frac{1}{\|p_{\tau_i+m-i}^{r,s}\|_2}Q_{\tau}^{r,s,m}(1-2x),
\end{equation*}
then the absolute convergence of \eqref{Jac-Sem} amounts to that of 
\begin{equation*}
\sum_{\tau_1 \geq \dots \geq \tau_m \geq 0}e^{-K_{\tau}^{r,s,m}(t/d)}\left[Q_{\tau}^{r,s,m}(1^m)\prod_{i=1}^m \frac{1}{\|p_{\tau_i+m-i}^{r,s}\|_2}\right]^2. 
\end{equation*}
By the virtue of the bound 
\begin{equation*}
V(\tilde{\tau}) \leq \prod_{i=1}^m[(\tau_i+m)(2\tau_i+2m+r+s+1)]^m
\end{equation*} 
and from the expression
\begin{equation*}
\|p_{\tau_i+m-i}^{r,s}\|_2^2 = \frac{1}{2(\tau_i+m-i)+r+s+1}\, \frac{\Gamma(r+\tau_i+m-i+1)\Gamma (s+\tau_i+m-i+1)}{\Gamma(\tau_i+m-i+1)\Gamma(\tau_i+m-i+1+r+s)},
\end{equation*}
it then suffices to prove the absolute convergence of the series 
\begin{multline*}
\sum_{\tau_1 \geq \dots \geq \tau_m \geq 0}e^{-K_{\tau}^{r,s,m}(t/d)} \Biggl(\prod_{i=1}^m[(\tau_i+m)(2\tau_i+2m+r+s+1)]^{2m} [2\tau_i+2m+r+s+1] \times \\ \frac{\Gamma(\tau_i+m-i+r+1)\Gamma(\tau_i+m-i+r+s+1)}{\Gamma(\tau_i+m-i+1)\Gamma(\tau_i+m-i+s+1)}\Biggr). 
\end{multline*}
Since this is a series of positive numbers, then we can bound it from above by the series over all the $m$-tuples $(\tau_1, \dots, \tau_m) \in \mathbb{N}^m$. Doing so leads to proving the absolute convergence of the series
\begin{align*}
\sum_{j \geq 0}e^{-j(j + r+s+1 + 2(m-i))(t/d)}&[(j+m)(2j+2m+r+s+1)]^{2m}\times \\& [2j+2m+r+s+1] (j+m-i+1)_r(j+m-i+s+1)_r,
 \end{align*}
 for any $1 \leq i \leq m$. But this holds true since 
\begin{equation*}
(j+m-i+1)_r(j+m-i+s+1)_r \sim (j+m-i+1)^r(j+m-i+s+1)^r, \quad j \rightarrow \infty.
\end{equation*}
From the mirror symmetry $q_n^{r,s}(-x) = (-1)^nq_n^{s,r}(x)$, it follows that $Q_{\tau}^{r,s,m}(-\phi) = (-1)^{|\tau|}Q_{\tau}^{r,s,m}(\phi)$ whence the absolute-convergence of the series \eqref{Jac-Sem} may be proved for $0 \leq r \leq s$ along the previous lines. As a matter of fact, if the hermitian matrix Jacobi process starts at the identity matrix $J_0 = {\it I}_m$, then Fubini Theorem yields 
\begin{align}\label{Integral}
\mathbb{E}(\mathrm{tr}[(J_{t/d})^n]) &= \int \left(\sum_{i=1}^m\lambda_i^n\right)G_t^{r,s,m}(1^m,\lambda)d\lambda \nonumber
\\& = \sum_{\tau_1 \geq \dots \geq \tau_m \geq 0}e^{-K_{\tau}^{r,s,m}(t/d)}P_{\tau}^{r,s,m}(1^m) \int \left(\sum_{i=1}^m\lambda_i^n\right)P_{\tau}^{r,s,m}(\lambda)W^{r,s,m}(\lambda)d\lambda.
\end{align}

\section{Proof of Theorem \ref{Obs}} In this section, we prove both Theorem \ref{Obs} and Corollary \ref{coll}.
The proof of the former relies mainly on the lemma below, where we determine the partitions having non zero contributions to the integral displayed in the right hand side of \eqref{Integral}.  
\subsection{Partitions} 
When $m=1$, $\tau$ is a nonnegative integer and $P_{\tau}^{r,s,1}$ reduces to the orthonormal one-dimensional Jacobi polynomial $P_{\tau}^{r,s}$ of degree $\tau$. In this case, the integral 
\begin{equation*}
\int_0^1 x^j P_{\tau}^{r,s}(x)x^r(1-x)^s dx 
\end{equation*}
vanishes unless $j \geq \tau$, since $x^j$ may be written as a linear combination of $P_{\tau}, \tau \leq j$. 

For general $m \geq 2$, the situation is quite similar. More precisely, fix $n < m$ and recall from \cite[page 68, exercise 10]{MD} the following expansion of the $n$-th power sum:
\begin{equation*}
 \sum_{i=1}^m\lambda_i^n=\sum_{k=0}^{n-1}(-1)^ks_{\alpha}(\lambda),
\end{equation*}
where 
\begin{equation*}
\alpha = \alpha(k,n) =  (n-k, 1^k):=(n-k, \underbrace{1,\ldots,1}_{k \mathrm{\,times}},\underbrace{0,\ldots,0}_{m-k-1 \mathrm{\,times}}), \quad 0 \leq k \leq n-1,
\end{equation*}  
are hooks of common weight 
\begin{equation*}
|\alpha| = \sum_{i=1}^m \alpha_i = n, 
\end{equation*}
and 
\begin{equation*}
s_{\alpha }(\lambda) = s_{\alpha}(\lambda_1,\ldots,\lambda_m)=\frac{\det(\lambda_j^{\alpha_i +m-i})_{i,j=1}^m}{\det(\lambda_j^{m-i})_{i,j=1}^m}
\end{equation*}
are the corresponding Schur polynomials. 

Recall also from \cite[page 37]{DG} the integral form of the Cauchy-Binet formula: for any probability measure $\kappa$ and any sequences $(\psi_i)_{i \geq 1}, (\phi_i)_{i \geq 1}$ of real-valued bounded functions,   
\begin{equation*}
\int \mathrm{det}(\psi_i(x_j))_{i,j =1}^m\mathrm{det}(\phi_i(x_j))_{i,j =1}^m \prod_{i=1}^m\kappa(dx_i) = m! \mathrm{det}\left(\int\psi_i(x)\phi_j(x)\kappa(dx)\right)_{i,j=1}^m.
\end{equation*}
We can now state the lemma alluded to above, where we use the ordering $\tau \subseteq \alpha$ meaning that $\tau_i \leq \alpha_i$ for all $1 \leq i \leq m$.
\begin{lem}\label{Partitions}
For any $k \leq n-1$, the integral 
\begin{equation*}
\int s_{\alpha}(\lambda)P_{\tau}^{r,s,m}(\lambda)W^{r,s,m}(\lambda)d\lambda 
\end{equation*}
vanishes unless $\tau \subseteq \alpha$.
\end{lem}

{\it Proof}: For sake of simplicity, let us omit in this proof the super-scripts and write simply $P_{\tau}, P_n, W$ instead of $P_{\tau}^{r,s,m}, P_n^{r,s}$, $W^{r,s,m}$ respectively. From the Cauchy-Binet formula, it follows that 
\begin{align*}
 \int s_{\alpha}(\lambda)P_{\tau}(\lambda)W(\lambda)d\lambda
 &= \frac{1}{m!} \int_{[0,1]^m} \mathrm{det}(\lambda_i^{\alpha_j + m-j}) \mathrm{det}(P_{\tau_j+m-j}(\lambda_i)) \prod_{i=1}^m\lambda_i^r(1-\lambda_i)^{s}d\lambda
\\& =  \mathrm{det}\biggl(\int_0^1x^{\alpha_j + m-j} P_{\tau_i+m-i}(x) x^r(1-x)^s dx\biggr)_{i,j=1}^m.
\end{align*}
Set
\begin{equation*}
A=(A_{ij})_{i,j=1}^m:=\biggl(\int_0^1x^{\alpha_j + m-j} P_{\tau_i+m-i}(x) x^r(1-x)^s dx\biggr)_{i,j=1}^m
\end{equation*}
and note that $\mathrm{det}(A) = 0$ if $\tau_m \geq 1$ since the last column is the null vector. Assuming $\tau_m=0, \tau_{m-1} \geq 1$ and expanding the determinant along the last column, then the same conclusion holds for the principal minor 
\begin{equation*}
(A_{ij})_{i,j=1}^{m-1}
\end{equation*}
and so on up to the principal minor of size $k+1$. Thus, $\det(A) = 0$ unless $\tau_i = 0$ for all $k+2 \leq i \leq m$. If $k=0$, then $A$ is a lower triangular matrix and $\det(A) = 0$ unless $\tau_1 \leq n$. Otherwise $1 \leq k \leq n-1$, and if $\tau_i \geq 2$ for some $2 \leq i \leq k+1$,  then $\tau_1 \geq \tau_2 \geq 2$ so that for any $j \geq 2$
\begin{equation*}
\tau_1 + m-1 \geq \tau_2+m-2 \geq m > \alpha_j + m-j.
\end{equation*}
From the orthogonality of the one-dimensional Jacobi polynomials, it follows that $A_{1j} = A_{2j} = 0$ for all $j \geq 2$ so that the first and the second row are proportional. Thus, $\det(A) = 0$ and we are left with the hooks 
\begin{equation*}
\tau = (\tau_1 \geq \underbrace{\tau_2 \geq\dots \geq \tau_{k+1}}_{\in \{0, 1\}} \geq \underbrace{0,\dots, 0}_{m-k-1 \mathrm{\,times}})
\end{equation*}
But if $\tau_1 > n-k \geq 1$ then the first row is the null vector and $\det(A) = 0$ as well. The lemma is proved. $\hfill \blacksquare$ 
\begin{rem}
We shall see below that the symmetric Jacobi polynomial has a `lower-triangular' expansion in the basis of Schur polynomials with respect to the ordering $\subseteq$. It is very likely that the inverse expansion of the Schur polynomial in the basis of symmetric Jacobi polynomials is also lower-triangular. In this case, the lemma would follow from the fact that symmetric Jacobi polynomials are mutually orthogonal with respect to $W^{r,s,m}$: 
\begin{equation*}
\int P_{\tau}^{r,s,m}(x)P_{\kappa}^{r,s,m}(x)W^{r,s,m}(x)dx  = 0
\end{equation*}
whenever the partitions $\tau$ and $\kappa$ are different.  
\end{rem}

Now we proceed to the end of the proof of Theorem \ref{Obs}. 

\subsection{Symmetric Jacobi polynomials associated with hooks}
Let $0 \leq k \leq n-1$ and $\tau \subseteq \alpha$ be a hook 
\begin{equation*}
\tau = (n-k-\delta, 1^{k-g}), \quad 0 \leq \delta \leq n-k-1, \quad 0 \leq g \leq k.
\end{equation*}
 For a partition $\mu$, we denote by 
\begin{equation*}
(z)_{\mu} = \prod_{i=1}^m (z-i+1)_{\mu_i} = \prod_{i=1}^m\frac{\Gamma(z-i+1+\mu_i)}{\Gamma(z-i+1)}
\end{equation*}
the generalized Pochhammer symbol.

From \cite{Las0} and \cite{Las1}, we dispose of an explicit expansion of $U_{\tau}^{r,s,m}$ in the Schur polynomial basis. More precisely,  by specializing \cite[Theorem 3]{Las1} to $\alpha=1$, we claim that
\begin{equation}\label{Expan}
U_{\tau}^{r,s,m}(\lambda) = \sum_{\mu \subseteq \tau}\frac{(-1)^{|\mu|}}{(r+m)_{\mu}}\binom{\tau}{\mu}C_{\mu}^{\tau}(r+s+2m)\frac{s_{\mu}(\lambda)}{s_{\mu}(1^m)}
\end{equation}
where if
\begin{equation*}
\mu = (n-k-\gamma, 1^{k-l}), \quad \delta \leq \gamma \leq n-k-1, \,\, g \leq l \leq k,
\end{equation*}
then
\begin{equation*}
\binom{\tau}{\mu} =  \binom{n-k-\delta-1}{\gamma-\delta}\binom{k-g}{l-g}\frac{(n-\delta - l)(n-g-\gamma) - (\gamma-\delta)(l-g)}{(n-\gamma-l)^2}
\end{equation*} 
is the generalized binomial coefficient (specialize \cite[Theorem 4]{Las0} to $\alpha=1$),  and where for any real $X$ (specialize  \cite[Theorem 6]{Las1} to $\alpha = 1$)
 \begin{multline}\label{Hook}
C_{\mu}^{\tau}(X) = \\ \left(X + \frac{(n-k-\delta)(n-k-\delta-1) - (k-g)(k-g+1)}{n-\delta - g}\right)\prod_{i=2}^{n-k-\gamma}(X+n-k-\delta+i-2)\prod_{i=1}^{k-l}(X-k+g-i).
\end{multline}

In order to prove Theorem \ref{Obs}, we need to compute
\[
\int s_{\alpha}(\lambda)U_{\tau}^{r,s,m}(\lambda)W^{r,s,m}(\lambda)d\lambda.
\]
With regard to \eqref{Jac-Sem1}, \eqref{Integral} and Lemma \ref{Partitions},

\begin{align*}
\int s_{\alpha}(\lambda)U_{\tau}^{r,s,m}(\lambda)&W^{r,s,m}(\lambda)d\lambda  =\sum_{\mu \subseteq \tau}\frac{(-1)^{|\mu|}}{(r+m)_{\mu}}\binom{\tau}{\mu}C_{\mu}^{\tau}(r+s+2m)  
\int \frac{s_{\alpha}(\lambda)s_{\mu}(\lambda)}{s_{\mu}(1^m)}W^{r,s,m}(\lambda)d\lambda 
\\& = \sum_{\mu \subseteq \tau}\frac{(-1)^{|\mu|}}{(r+m)_{\mu}s_{\mu}(1^m)}\binom{\tau}{\mu}C_{\mu}^{\tau}(r+s+2m)  
 \det\left(\int_0^1x^{\alpha_i + \mu_j+2m - i- j +r}(1-x)^sdx\right)_{i,j=1}^m
\\& = \sum_{\mu \subseteq \tau}\frac{(-1)^{|\mu|}}{(r+m)_{\mu}s_{\mu}(1^m)}\binom{\tau}{\mu}C_{\mu}^{\tau}(r+s+2m)  
 \det\left(\beta(\alpha_i + \mu_j+2m - i- j +r +1,s+1)\right)_{i,j=1}^m.
\end{align*}
The formula displayed in Theorem \ref{Obs} follows after setting 
\begin{eqnarray}\label{coeff1}
a_{\tau}^{r,s,m} &:=& \left\{V(\tilde{\tau}) \prod_{1 \leq i < j \leq m} \frac{1}{(r+j-i)i} \prod_{i=1}^m \frac{ p_{\tau_i+m-i}^{r,s}(0)}{||p_{\tau_i+m-i}^{r,s}||_2}\right\}^2, \\  
\label{coeff2} 
b_{\mu,\tau}^{r,s,m} &:=& \frac{(-1)^{|\mu|}}{(r+m)_{\mu}s_{\mu}(1^m)}\binom{\tau}{\mu}C_{\mu}^{\tau}(r+s+2m).
\end{eqnarray}

\begin{rem}
The product $s_{\alpha}s_{\mu}$ is linearized via the Littlewood-Richardson coefficients (\cite{MD}, p.142) as: 
\begin{equation*}
s_{\alpha}(\lambda)s_{\mu}(\lambda) = \sum_{\kappa}c_{\alpha \mu}^{\kappa}s_{\kappa}(\lambda),
\end{equation*}
where the summation is over the set of partitions $\{\kappa \supseteq \alpha, \kappa \supseteq \mu, |\alpha| + |\mu| = |\kappa|\}$. Thus 
\begin{equation}\label{ProdSchur}
\int s_{\alpha}(\lambda)s_{\mu}(\lambda)W^{r,s,m}(\lambda)d\lambda =  \sum_{\kappa}c_{\alpha \mu}^{\kappa} \int s_{\kappa}(\lambda)W^{r,s,m}(\lambda)d\lambda
\end{equation}
and the value of the integral in the right hand side is an instance of Kadell's integral (see Exercice 7, p.385 in \cite{MD}):
\begin{equation*}
\int s_{\kappa}(\lambda)W^{r,s,m}(\lambda)d\lambda =  \prod_{1\leq i < j \leq m}(\kappa_i-\kappa_j+j-i)\prod_{i=1}^m\frac{\Gamma(\kappa_i+r+m-i+1)\Gamma(s+m-i+1)}{\Gamma(\kappa_i+r+s+2m-i+1)}.
\end{equation*}
However, up to our best knowledge, there is no simple formula for $c_{\alpha \mu}^{\kappa}$ except when $\mu$ is a partition with one row or one column\footnote{This is referred to as Pieri formula.}. For that reason, we preferred the use of the Cauchy-Binet formula when evaluating \eqref{ProdSchur}. Nonetheless, if $\tau = (0)$ is the null partition then $\mu = (0)$ and the left-hand side of \eqref{ProdSchur} reduces to Kadell's integral. Moreover, $b_{0,0}^{r,s,m} = 1, U_{(0)}^{r,s,m} = 1$ and $a_{(0)}^{r,s,m}$ is exactly the normalizing constant of $W^{r,s,m}$ whose multiplicative inverse is a special instance of the value of the Selberg integral (see e.g. \cite{CDLV}). Consequently, if we let $t \rightarrow \infty$ in \eqref{Obs}, then the only non-vanishing term corresponds to 
$\tau = (0)$ and as such, we retrieve the moments of $W^{r,s,m}$ (which is the stationary distribution of the eigenvalues process $(\lambda(t))_{t \geq 0}$) derived in \cite{CDLV}.   
\end{rem}

\subsection{The case $s=0$: proof of Corollary \ref{coll}}
Specializing Theorem \ref{Obs} with $s=0$, then the Cauchy-determinant yields:
\begin{align*}
\det\left(\beta(\alpha_i + \mu_j+2m - i- j +r +1,1)\right)_{i,j=1}^m  &= \det\left(\frac{1}{\alpha_i + \mu_j+2m - i- j +r +1}\right)_{i,j=1}^m
\\& = \frac{\prod_{1\leq i < j \leq m}(\alpha_i - \alpha_j+j-i)(\mu_i-\mu_j+j-i)}{\prod_{i,j=1}^m(\alpha_i+\mu_j+2m-i-j+r+1)}.
\end{align*}
Besides, the Weyl dimension formula
\begin{equation*}
s_{\mu}(1^m) =  \prod_{1\leq i < j \leq m}\frac{(\mu_i-\mu_j+j-i)}{j-i}
\end{equation*}
and the equality 
\begin{equation*}
\prod_{1 \leq i < j \leq m}[(r+j-i)i] = \prod_{i=1}^m\frac{\Gamma(r+m-i+1)}{\Gamma(r+1)}\prod_{1 \leq i < j \leq m}(j-i)
\end{equation*}
entail
\begin{equation*}
\prod_{1\leq i < j \leq m}\frac{(\alpha_i - \alpha_j+j-i)(\mu_i-\mu_j+j-i)}{[(r+j-i)i]^2}  = \prod_{i=1}^m\left\{\frac{\Gamma(r+1)}{\Gamma(r+m-i+1)}\right\}^2s_{\alpha}(1^m)s_{\mu}(1^m).
\end{equation*}
Formula \eqref{Formula1} in Corollary \ref{coll} follows then from the equality
\begin{equation*}
\prod_{1 \leq i < j \leq m}(j-i) = \prod_{i=1}^m\Gamma(m-i+1)
\end{equation*}
together with
\begin{equation*}
p_{\tau_i+m-i}^{r,0}(0) = \frac{\Gamma(r+1+\tau_i+m-i)}{\Gamma(r+1)\Gamma(\tau_i+m-i+1)}, \quad ||p_{\tau_i+m-i}^{r,0}||^2_2 = \frac{1}{2(\tau_i+m-i)+r+1}.
\end{equation*}
The second formula in the corollary is obvious.

\section{Asymptotics}
The purpose of this section is to determine the limits of various terms appearing in \eqref{Formula1} under the assumption that the limits \eqref{Limits} exist. Doing so is the crucial step in our future investigations aiming in particular to derive the moments \eqref{Mom} as limits of their matrix analogues and more generally to derive an expression for $\mathcal M_n(t, \eta, \theta)$. We start with 
\begin{equation*}
\lim_{m \rightarrow \infty} \frac{1}{d(m)}K_{\tau}^{r(m),s(m),m} = \lim_{m \rightarrow \infty} \frac{1}{d(m)} \sum_{i=1}^{l(\tau)}\tau_i(\tau_i+d(m)+1-2i) = |\tau|
\end{equation*} 
which holds for any hook $\tau$ of weight $|\tau| \leq n$. Next, we prove the following lemma:
\begin{lem} 
Let $\tau$ be a hook of weight $|\tau| \leq n$ and let $\mu \subseteq \tau$. Then 
\begin{equation*}
\lim_{m \rightarrow \infty} b_{\mu,\tau}^{r(m),s(m),m}s_{\mu}(1^m) = \frac{(-1)^{|\mu|}}{\theta^{|\mu|}}\binom{\tau}{\mu}. 
\end{equation*}
In particular, 
\begin{equation*}
\lim_{m \rightarrow \infty}U_{\tau}^{r(m),s(m),m}(1^m) = \left(1-\frac{1}{\theta}\right)^{|\tau|}.
\end{equation*}
\end{lem}
{\it Proof}: Since $l(\mu) \leq n < m$, then the generalized Pochammer symbol splits as
\begin{align*}
(r(m)+m)_{\mu} = \prod_{i=1}^{l(\mu)}(p(m)-i+1)_{\mu_i} = (p(m))_{\mu_1}\prod_{i=2}^{l(\mu)}(p(m)-i+1).
\end{align*}
Thus $(r(m)+m)_{\mu} \sim p(m)^{|\mu|}$ as $m \rightarrow \infty$. On the other hand, it is obvious from \eqref{Hook} that 
\begin{equation*}
C_{\mu}^{\tau}(r(m)+s(m)+2m) = C_{\mu}^{\tau}(d(m)) \sim d(m)^{|\mu|} \mathrm{\ as\ }\, m \rightarrow \infty. 
\end{equation*}
Hence, we get from \eqref{coeff2}: 
\begin{equation*}
\lim_{m \rightarrow \infty} b_{\mu,\tau}^{r(m),s(m),m}s_{\mu}(1^m) =  \lim_{m \rightarrow \infty} \frac{(-1)^{|\mu|}}{(r(m)+m)_{\mu}}\binom{\tau}{\mu}C_{\mu}^{\tau}(r(m)+s(m)+2m),
= \frac{(-1)^{|\mu|}}{\theta^{|\mu|}}\binom{\tau}{\mu}
\end{equation*}
and from \eqref{Expan}:
\begin{equation*}
\lim_{m \rightarrow \infty}U_{\tau}^{r(m),s(m),m}(1^m) = \sum_{\mu \subseteq \tau} (-1)^{|\mu|}\binom{\tau}{\mu} \left(\frac{1}{\theta}\right)^{|\mu|} = \frac{s_{\tau}(1-(1/\theta), \dots, 1-(1/\theta))}{s_{\tau}(1^{l(\tau)})},
\end{equation*}
where the last equality follows from the generalized binomial Theorem (\cite{Lass}). The lemma follows from the homogeneity of the Schur polynomials. $\hfill \blacksquare$ 

Now, assume $s(m) = 0$ and note that this assumption yields in the large $m$-limit the relation 
\begin{equation*}
\theta(1+\eta) = 1 \quad \Leftrightarrow \quad \eta = \frac{1-\theta}{\theta}.
\end{equation*}
If $l(\mu) \leq l(\tau) \leq l(\alpha)\leq n<m$ are the lengths of the partitions $\mu \subseteq \tau \subseteq \alpha$ respectively, then the following cancellations occur: 
\begin{equation*}
\prod_{i=1}^m \frac{\Gamma(r+\tau_i+m-i+1)\Gamma(m-i+1)}{\Gamma(\tau_i+m-i+1)\Gamma(r+m-i+1)} = \prod_{i=1}^{l(\tau)}\frac{\Gamma(r+\tau_i+m-i+1)\Gamma(m-i+1)}{\Gamma(\tau_i+m-i+1)\Gamma(r+m-i+1)},
\end{equation*}
\begin{equation*}
\prod_{i=1}^m\frac{[2(\tau_i+m-i)+r+1]}{(\alpha_i+\mu_i +2m-2i+r+1)} = \prod_{i=1}^{l(\alpha)}\frac{[2(\tau_i+m-i)+r+1]}{(\alpha_i+\mu_i +2m-2i+r+1)}, 
\end{equation*}
and 
\begin{align*}
\frac{\prod_{l(\alpha)+1 \leq i < j \leq m}(\tau_i+\tau_j+2m-i-j+r+1)^2}{\prod_{l(\alpha)+1\leq i \neq j\leq m}(\alpha_i+\mu_j+2m-i-j+r+1)} =  
\prod_{l(\alpha)+1 \leq i \neq j \leq m} \frac{(\tau_i+\tau_j+ 2m-i-j+r+1)}{(\alpha_i+\mu_j+2m-i-j+r+1)} =1.
\end{align*}
As a result, 
\begin{align*}
\lim_{m \rightarrow \infty} \prod_{i=1}^m \frac{\Gamma(r(m)+\tau_i+m-i+1)\Gamma(m-i+1)}{\Gamma(\tau_i+m-i+1)\Gamma(r(m)+m-i+1)} & = \prod_{i=1}^{l(\tau)} \lim_{m \rightarrow \infty} \left(\frac{p(m)}{m}\right)^{\tau_i} = \frac{1}{\eta^{|\tau|}}
= \Bigr( \frac{\theta}{1-\theta}\Bigl)^{|\tau|},
\end{align*}
and similarly

\begin{align*}
 \lim_{m \rightarrow \infty}\prod_{i=1}^m\frac{[2(\tau_i+m-i)+r(m)+1]}{(\alpha_i+\mu_i +2m-2i+r(m)+1)} = 1,
 \end{align*}

 \begin{equation*}
\lim_{m \rightarrow \infty} \frac{\prod_{1\leq i < j \leq l(\alpha)}(\tau_i+\tau_j+2m-i-j+r(m)+1)^2}{\prod_{1\leq i \neq j\leq l(\alpha)}(\alpha_i+\mu_j+2m-i-j+r(m)+1)} = 
\lim_{m \rightarrow \infty} \prod_{1\leq i \neq j \leq l(\alpha)}\frac{(\tau_i+\tau_j+2m-i-j+r(m)+1)}{(\alpha_i+\mu_j+2m-i-j+r(m)+1)} = 1.
 \end{equation*}
Finally, consider the product
\begin{equation*}
\frac{\displaystyle \prod_{\substack{1\leq i \leq l(\alpha)\\ l(\alpha)+1\leq j \leq m}}(\tau_i+2m-i-j+r(m)+1)^2}{\displaystyle \prod_{\substack{1\leq i \leq l(\alpha) \\ l(\alpha)+1\leq j \leq m}}(\alpha_j + \mu_i +2m-i-j+r(m)+1)
\prod_{\substack{1\leq j \leq l(\alpha)\\ l(\alpha)+1\leq i \leq m}}(\alpha_i + \mu_j +2m-i-j+r(m)+1)}.
\end{equation*}
It can be rewritten as
\begin{equation*}
\prod_{1\leq i \leq l(\alpha)}\prod_{l(\alpha)+1\leq j \leq m}\frac{(\tau_i+2m-i-j+r(m)+1)^2}{(\alpha_i+2m-i-j+r(m)+1)(\mu_i+2m-i-j+r(m)+1)}
\end{equation*}
which shows that it is equivalent to $[d(m)]^{2|\tau|-|\alpha| - |\mu|}$ as $m \rightarrow \infty$. Indeed, recall $r(m)+2m = p(m)+m = d(m)$ and consider
\begin{equation*}
\prod_{1\leq i \leq l(\alpha)}\prod_{l(\alpha)+1\leq j \leq m}\frac{(\tau_i+d(m)-i-j+1)}{(d(m)-i-j+1)} = \prod_{1\leq i \leq l(\tau)}\prod_{l(\alpha)+1\leq j \leq m}\frac{(\tau_i+d(m)-i-j+1)}{(d(m)-i-j+1)}. 
\end{equation*}
Then the terms corresponding to $i=1$ are
\begin{equation*}
\frac{(d(m)-l(\alpha)+\tau_1-1)(d(m)-l(\alpha)+\tau_1-2)\cdots(d(m)-l(\alpha))(d(m)-l(\alpha)-1)\cdots (d(m)-m)}{(d(m)-l(\alpha)-1)\cdots(d(m)-m)}
\end{equation*}
which reduces to 
\begin{equation*}
\prod_{j=0}^{\tau_1 - 1} (d(m)-l(\alpha)+j) \sim d(m)^{\tau_1}, \quad m \rightarrow \infty.
\end{equation*}
Consequently 
 \begin{equation*}
\prod_{1\leq i \leq l(\alpha)}\prod_{l(\alpha)+1\leq j \leq m}\frac{(\tau_i+d(m)-i-j+1)}{(d(m)-i-j+1)} \sim d(m)^{|\tau|}, \quad m \rightarrow \infty.
\end{equation*}
The same reasoning shows that 
\begin{equation*}
\prod_{1\leq i \leq l(\alpha)}\prod_{l(\alpha)+1\leq j \leq m}\frac{(\alpha_i+d(m)-i-j+1)}{(d(m)-i-j+1)} \sim d(m)^{|\alpha|}, \quad m \rightarrow \infty,
\end{equation*}

\begin{equation*}
\prod_{1\leq i \leq l(\alpha)}\prod_{l(\alpha)+1\leq j \leq m}\frac{(\mu_i+d(m)-i-j+1)}{(d(m)-i-j+1)} \sim d(m)^{|\mu|}, \quad m \rightarrow \infty,
\end{equation*}
whence the claimed equivalence follows.

Summing up, all the terms of the finite sum in the right hand side of formula \eqref{Formula1} admit finite limits except $s_{\alpha}(1^m)$ and $s_{\mu}(1)^m$. Since the latter are equivalent to $d(m)^{|\alpha|}$ and to $d(m)^{|\mu|}$ respectively as $m \rightarrow \infty$ and due to the presence of alternating signs, taking the limit as $m \rightarrow \infty$ in formula \eqref{Formula1} leads to an indeterminate limit. 
To solve this problem, one needs to seek some cancellations in a similar fashion this was done for the unitary Brownian motion (\cite{Biane}). 

\end{document}